\documentclass[11pt,twoside,reqno]{amsart}
\usepackage{color}
\usepackage{amsfonts}
\usepackage{}

\usepackage[OT1]{fontenc}
\usepackage{type1cm}
\usepackage{amssymb}
\usepackage{cite}
\usepackage{xcolor}
\usepackage[left=2.3cm,top=3cm,right=2.3cm]{geometry}

\numberwithin{equation}{section}

\theoremstyle{plain}

\newtheorem{thm}{Theorem}

\newtheorem{proposition}{Proposition}

\newtheorem{lemma}{Lemma}

\theoremstyle{remark}
\newtheorem{remark}{Remark}

\theoremstyle{definition}

\newtheorem{defn}{Definition}

\begin{document}

\title{Assouad type dimensions of  homogeneous Moran sets and Cantor-like sets}

\author{Yanzhe Li$^{*}$}
\address[Yanzhe Li]{College of Mathematics and Information Science, Guangxi University, Nanning, 530004, P.~R. China}
\email{lyzkbm@163.com}

\author{Jun Li}
\address[Jun Li]{College of Mathematics and Information Science, Guangxi University, Nanning, 530004, P.~R. China}
\email{836164679@qq.com}

\author{Shuang Liang}
\address[Shuang Liang]{College of Mathematics and Information Science, Guangxi University, Nanning, 530004, P.~R. China}
\email{liangshuangkw@163.com}

\author{Manli Lou}
\address[Manli Lou]{Department of Basic Courses, Guangzhou Maritime University, Guangzhou, 510725, P.~R. China}
\email{loumanli@126.com}

\thanks{$*$ Corresponding author. }
\thanks{This work is supported by National Natural Science Foundation of China (No. 12461015), Guangxi Natural Science Foundation (2020GXNSFAA297040)and Guangdong Natural Science Foundation (2018A030313971).}


\begin{abstract}
In this paper, we give the Assouad dimension formula and the upper bound of the lower dimension for homogeneous Moran sets under the condition $\sup_{k\ge 1}\{n_{k}\}<+\infty$. We also give the Assouad spectrum and the lower spectrum formulas for Cantor-like sets.
\end{abstract}
\maketitle

\section{Introduction}
Let $(X, d)$ be a metric space, for any set $F\subset (X, d)$, define the Assouad dimension of $F$(denoted by $\dim_{A}F$) as
\begin{equation}
\begin{split}
\dim_{A}F=\inf \bigg\{s\geq 0:\ & \text{there exist} \ b,\ c>0 \ \text{such that for any}\ \ 0<r< R< b, \ x\in F, \\& N_{r}\big(B(x,R)\cap F\big)\leq c\big(\frac{R}{r}\big)^{s} \bigg\},   \label{Assouad1}
\end{split}
\end{equation}
where $N_{r}(E)$  denotes the smallest number of the balls with radius $r$ which can cover the set $E\subset X$.
The definition of the Assouad dimension  was first introduced by Assouad \cite{As1} \cite{As2} and used for the study of the embedding
theories in the metric space. Unlike the Hausdorff, box and packing dimensions, the Assouad dimension  mainly describes the local strcture of the sets(the "density" of the most "dense" part of the sets), it plays effective roles on the study of fractal geometry, more details for the properties and the applications of the Assouad dimension can be found in \cite{Fra21,Fra14,FT18,Luu}.

 For any subset $F$ in $\mathbb{R}^{n}$, we have $\dim_{H}F \leq \dim_{P}F \leq \dim_{A}F$, and we also have $\overline{\operatorname{dim}}_{B}F \leq \dim_{A}F$ for any bounded set $F$, see in \cite{Fra21}. Furthermore, by  \cite{Luu}, we have for any Ahlfors$-$David $s-$regular set $F$, the Assouad dimension of $F$ is $s$, which is coincided with the Hausdorff dimension  and  the packing dimension of $F$. Calculate the Assouad dimensions of the sets which are not Ahlfors$-$David regular is  always a difficult work. Fraser\cite{Fra14} and Mackay \cite{Ma11}   discussed the Assouad dimensions for some  quasi$-$self$-$similar and self$-$affine sets. Li, Li, Miao and Xi \cite{LLMX} obtained some results of  the Assouad dimensions for Cantor-like sets and some Moran sets. Chen \cite{C19}, Peng, Wang and Wen \cite{PWW17},  Yang and Du \cite{YD20}, Dai, Dong and Wei \cite{DDW23} discussed the Assouad dimensions for some fractal sets in $\mathbb{R}$, such as the Moran cut-out sets under some conditions and some special homogeneous Moran sets. In this paper, we study the Assouad dimensions for the homogeneous Moran sets which satisfy  $\sup_{k\ge 1}\{n_{k}\}<+\infty$, and obtain the Assouad dimension formula, which generalizes the results in \cite{LLMX}, \cite{PWW17}, \cite{YD20} and \cite{DDW23}.
 \bigskip

As a dual to the Assouad dimension, Larman \cite{Lar1} introduced the lower dimension. Let $(X, d)$ be a metric space, for any set $F\subset (X, d)$, define the lower dimension of $F$(denoted by $\dim_{L}F$) as
\begin{equation}
\begin{split}
\dim_{L}F=\sup \bigg\{s\geq 0:\ & \text{there exist} \ b,\ c>0 \ \text{such that for any}\ \ 0<r< R< b, \ x\in F, \\& N_{r}\big(B(x,R)\cap F\big)\geq c\big(\frac{R}{r}\big)^{s} \bigg\}.   \label{Assouad2}
\end{split}
\end{equation}
Similar to the  Assouad dimension, the lower dimension  mainly describes the local strcture of the sets(the "sparsity" of the most "sparse" part of the sets), it also plays effective roles on the study of fractal geometry. But unlike the Assouad dimension, the lower dimension  does not have the monotone property and the finitely stable property,  more details for the properties and the applications of the lower dimension can be found in  \cite{Fra21}.

By \cite{Lar1}, we have $\dim_{L}F \leq \underline{\dim}_{B}F\leq \overline{\dim}_{B}F \leq \dim_{A}F$ for any  totally bounded set $F$. In general, since the lower dimension  does not have the monotone property and the finitely stable property, it is alway difficult to do the study of the lower dimensions for fractal sets. Chen, Wu, Wei\cite{CWW17} discussed the lower dimensions for some Moran sets. Chen\cite{C19} , Yang, Li and Hu\cite{YL22} obtained some results of the lower dimension formulas for some Moran cut-out sets and some homogeneous Moran sets. In this paper, we study the lower dimensions for the homogeneous Moran sets which satisfy  $\sup_{k\ge 1}\{n_{k}\}<+\infty$, and obtain an upper bound of the lower dimension.

\bigskip
By the definitions above, the Assouad dimension and the lower dimension are infimum and supremum depending on the two independent scales $r$ and $R$, but there is no information that how the scales approach the values. To see how the infimum and supremum depend on the scales, Fraser and Han \cite{FH18}  restricted  the relationship
between  $r$ and $R$  with $\frac{\log R}{\log r}=\theta(0<\theta <1)$  and gived the definitions of the Assouad spectrum and the lower spectrum. Let $(X, d)$ be a metric space, for any set $F\subset (X, d)$ and any $0<\theta <1$,  define the Assouad spectrum  of $F$(denoted by $\dim_{A}^{\theta}F$) as
\begin{equation}
\begin{split}
 \dim_{A}^{\theta}F=\inf \bigg\{s\geq 0:\ & \text{there exist} \ b,\ c>0 \ \text{such that for any}\ \ 0< R< b,\ x\in F, \\& N_{R^{\frac{1}{\theta}}}\big(B(x,R)\cap F \big)\leq c\big(\frac{R}{R^{\frac{1}{\theta}}}\big)^{s} \bigg\},    \label{Assouads1}
\end{split}
\end{equation}
and define the lower spectrum  of $F$(denoted by $\dim_{L}^{\theta}F$) as
\begin{equation}
\begin{split}
 \dim_{L}^{\theta}F=\sup \bigg\{s\geq 0:\ & \text{there exist} \ b,\ c>0 \ \text{such that for any}\ \ 0< R< b,\ x\in F, \\& N\big(B(x,R)\cap F, R^{\frac{1}{\theta}}\big)\geq c\big(\frac{R}{R^{\frac{1}{\theta}}}\big)^{s} \bigg\}.   \label{lowers1}
\end{split}
\end{equation}
If $F$ is a totally bounded set, it is easy to see that $ \dim_{A}^{\theta}F$ and  $\dim_{L}^{\theta}F$ are continuous functions for any $0<\theta <1$, and we have
\begin{equation*}
\overline{\dim}_{B}F \leq \dim_{A}^{\theta}F \leq \min\bigg\{\frac{\overline{\dim}_{B}F}{1-\theta}, \dim_{A}F \bigg\}
\end{equation*}
and
\begin{equation*}
\dim_{L}F\leq \dim_{L}^{\theta}F\leq \underline{\dim}_{B}F.
\end{equation*}
More details for the properties and the applications of the Assouad spectrum and the lower spectrum can be found in  \cite{Fra21} and \cite{FH18}.

Fraser \cite{Fra14}  studied the Assouad spectrums and the lower spectrums for some classical self$-$similar sets and self$-$conformal sets. Chen \cite{C19} obtained some results of the Assouad spectrums and the lower spectrums for some special Moran cut-out sets. Yang and Du \cite{YD20} obtained  some results of the  Assouad spectrums,   Yang, Li and Hu \cite{YL22} obtained  some results of the  lower spectrums for a class of special homogeneous Moran sets which is called homogenoeous perfect sets under some conditions. Chen, Wu, Chang\cite{CWC} studied some properties of the lower spectrums for uniformly perfect sets. In this paper, we study the Assouad spectrums and the lower spectrums for the Cantor-like sets, and obtain the Assouad spectrum and the lower spectrum formulas for the Cantor-like sets, which generalize the results in \cite{YD20} and \cite{YL22}.
\bigskip

The paper is organized as follows. In Section 2, we recall the definitions of homogeneous Moran sets and the Cantor-like sets. Our main results are stated in Section 3. Section 4 and Section 5 are the proofs of our main results.

\bigskip
\bigskip
\section{Preliminaries}
\label{sec:pre1}
We recall the definition of the homogeneous Moran sets.
Let  the sequences $\{n_{k}\}_{k\geq 1}\subset \mathbb{Z}^+$ and $%
\{c_{k}\}_{k\geq 1}\subset \mathbb{R}^+$
with $n_{k}\geq 2$, $n_{k}c_{k}<1$ for any $k\geq 1$%
. Set $\Omega _{0}=\{\emptyset\}$, $\Omega _{k}=\{\sigma_{1}\sigma_{2}\cdots \sigma_{k}: 1\leq \sigma_{j} \leq n_{j}, 1\leq j \leq k\}$ for any $k\geq 1$, and define $\Omega=\bigcup_{k=0}^{\infty}\Omega_{k}$. For any $k\geq 1$, $m\geq 1$, if $\sigma=\sigma_{1}\sigma_{2}\cdots \sigma_{k} \in \Omega _{k}$,$\tau=\tau_{1}\tau_{2}\cdots \tau_{m}(1 \leq \tau_{j}\leq n_{k+j}, 1 \leq j \leq m)$ , then write $\sigma\ast \tau=\sigma_{1}\cdots \sigma_{k}\tau_{1}\cdots\tau_{m} \in \Omega _{k+m}$.
\begin{defn} (Homogeneous Moran sets \cite{Wen00})  \label{hps}

Let $I\neq \emptyset$ be a closed interval in $\mathbb{R}$, we say the collection of the closed intervals $\mathcal{F}=\{I_{\sigma}\subset I:\sigma\in \Omega\}$ fulfills the homogeneous Moran structure if it satisfies:

\begin{enumerate}

\item[\textup{(i)}] $I_{\emptyset}=I$;
\item[\textup{(ii)}]for any $k\geq 1$ and $\sigma \in \Omega _{k-1}$, $ I_{\sigma \ast 1\text{,}}$\ $ I_{\sigma \ast 2\text{,}}\cdots $, $I_{\sigma \ast n_{k}}$ are closed intervals with $\bigcup_{i=1}^{ n_{k}}I_{\sigma \ast i}\subset I_{\sigma}$ and $\text{int}(I_{\sigma \ast i_{1}})\cap \text{int}(I_{\sigma \ast i_{2}})=\emptyset$ for any $1\le i_{1}<i_{2}\le n_{k}$;
where $\text{int}(A)$ denotes the interiors of the set $A\subset \mathbb{R}$;
\item[\textup{(iii)}] for any $k\geq 1$, $\sigma \in \Omega _{k-1}$, and $1\leq i\leq n_{k}$,
\begin{equation*}
\frac{|I_{\sigma\ast i}|}{|I_{\sigma }|}=c_{k},
\end{equation*}
where $|A|$ denotes the diameter of the set $A\subset \mathbb{R}$.

\end{enumerate}

If the collection $\mathcal{F}=\{I_{\sigma}:\sigma\in \Omega\}$ fulfills the homogeneous Moran structure, then we call
\begin{equation}
E=E(\mathcal{F}):=\bigcap_{k \geq 1}\bigcup_{\sigma\in\Omega_{k}}I_{\sigma}
\end{equation}
a homogeneous Moran set defined by $\mathcal{F}$. For any $k\ge 1$, if $\sigma\in \Omega_{k}$, then $I_{\sigma}\in \mathcal{F}$ is called a $k-$level basic interval of $E$. Denote by $\mathcal{M}(I, \{n_{k}\}, \{c_{k}\})$ the class of all homogeneous Moran sets which are associated with $I, \{n_{k}\}$ and $ \{c_{k}\}$.
\end{defn}
\begin{remark}
Without loss of generality, we  assume that $I=[0,1]$.
\end{remark}

Let
\begin{equation*}
N_{k}=\prod_{i=1}^{k}n_{i},\ \ \ \  \delta_{0}=1,\ \ \ \  \delta_{k}=\prod_{i=1}^{k}c_{i}
\end{equation*}
for any $k\ge 1$, then $N_{k}$ is  the total number of the $k-$level basic intervals of $E$ and  $\delta_{k}$ is   the length of any $k-$level basic interval of $E$.

\begin{remark}
Any uniform Cantor set(see in \cite{PWW17}) is a homogeneous Moran set, any homogeneous perfect set(see in \cite{YD20}) is a homogeneous Moran set, any parameterized  homogeneous Moran set(see in \cite{DDW23}) is a homogeneous Moran set .
\end{remark}

More details for the homogeneous Moran sets can be found in \cite{Wen00,FWW97,KLS,LW11}.
\bigskip
\bigskip

The definition of the Cantor-like sets is introduced by  Li, Li, Miao and Xi in  \cite{LLMX} as follows.
\begin{defn} (Cantor-like sets \cite{LLMX})  \label{CL}

Let  the sequence $\{a_{k}\}_{k\geq 1}\subset \mathbb{R}^+$ satisfying $\sum_{k=1}^{\infty}
a_{k}<+\infty$. Let $\{n_{k}\}_{k\geq 1}\subset \mathbb{Z}^+$ and $%
\{c_{k}\}_{k\geq 1}\subset \mathbb{R}^+$ be the sequences
with $n_{k}\geq 2$, $0<c_{k}<1$ and $c_{\ast}=\inf_{k}\{c_{k}\}>0$ for any $k\geq 1$. Set $\Omega _{0}=\{\emptyset\}$, $\Omega _{k}=\{\sigma_{1}\sigma_{2}\cdots \sigma_{k}: 1\leq \sigma_{j} \leq n_{j}, 1\leq j \leq k\}$ for any $k\geq 1$, and define $\Omega=\bigcup_{k=0}^{\infty}\Omega_{k}$. For any $k\geq 1$, $m\geq 1$, if $\sigma=\sigma_{1}\sigma_{2}\cdots \sigma_{k} \in \Omega _{k}$,$\tau=\tau_{1}\tau_{2}\cdots \tau_{m}(1 \leq \tau_{j}\leq n_{k+j}, 1 \leq j \leq m)$ , then write $\sigma\ast \tau=\sigma_{1}\cdots \sigma_{k}\tau_{1}\cdots\tau_{m} \in \Omega _{k+m}$.

Let $J\subset \mathbb{R}^{n}(n\ge 1)$ be a closed set  with $\text{int}(J) \neq \emptyset$. For any $k\ge 1$ and $\sigma\in \Omega_{k-1}$, if
$$J_{\sigma*1},J_{\sigma*2},\cdots,J_{\sigma*n_{k}}\subset J_{\sigma}$$
are geometrically similar to $J_{\sigma}$ with
$$c_{k}(1-a_{k})\le\frac{|J_{\sigma*j}|}{|J_{\sigma}|}\le c_{k}(1+a_{k})$$
for any $1\le j\le n_{k}$, and $\text{int}(J_{\sigma \ast i_{1}})\cap \text{int}(J_{\sigma \ast i_{2}})=\emptyset$ for any $1\le i_{1}<i_{2}\le n_{k}$, then
\begin{equation}
E=\bigcap_{k \geq 1}\bigcup_{\sigma\in\Omega_{k}}J_{\sigma}
\end{equation}
is called a Cantor-like set. Denote by $\mathcal{C}(J,  \{c_{k}\},\{n_{k}\},\{a_{k}\})$ the class of all Cantor-like sets  associated with $J, \{n_{k}\}, \{c_{k}\}$ and $\{a_{k}\}$.
\end{defn}

For any $k \ge 1$, we also write
\begin{equation*}
N_{k}=\prod_{i=1}^{k}n_{i},\ \ \ \  \delta_{0}=1,\ \ \ \  \delta_{k}=\prod_{i=1}^{k}c_{i}.
\end{equation*}

\begin{remark}
If $E\in \mathcal{M}(I, \{n_{k}\}, \{c_{k}\})$ is a homogeneous Moran set with $c_{\ast}=\inf_{k}\{c_{k}\}>0$, then $E$ is a Cantor-like set with $E\in \mathcal{C}(I,  \{c_{k}\},\{n_{k}\},\{a_{k}\})$ and $a_{k}=0$ for all $k\ge 1$. Cantor-like set may not be Moran set(see in \cite{Wen00}).
\end{remark}
\bigskip
\bigskip

\section{Main results}
 We state our main results as follows.

\begin{thm}\label{thm1}
Let $E\in \mathcal{M}(I, \{n_{k}\}, \{c_{k}\})$, if $\ \sup_{k\ge 1}\{n_{k}\}<+\infty$ , then
\begin{equation*}
\dim_{A}E=\limsup\limits_{l\rightarrow +\infty}\sup\limits_{k\geq 1}\frac{\log n_{k+1}\cdots n_{k+l}}
{-\log c_{k+1}\cdots c_{k+l}}.
\end{equation*}
\end{thm}

\begin{remark}

Any homogeneous Cantor set $E \in \mathcal{C}(\{n_{k}\}, \{c_{k}\})$, any homogeneous perfect set $F \in \mathcal{J}(I, \{n_{k}\}, \{c_{k}\}, \{\eta_{k,j}\} )$ and any parameterized  homogeneous Moran set $G \in \mathcal{M}(I, \{n_{k}\}, \{c_{k}\}, \{\varepsilon_{k}(t)\} )$ are  homogeneous Moran sets, where $E, F, G\in \mathcal{M}(I, \{n_{k}\}, \{c_{k}\})$. If $c_{\ast}=\inf_{k}\{c_{k}\}>0$, then $\ \sup_{k\ge 1}\{n_{k}\}<+\infty$. Thus Theorem \ref{thm1} generalizes Theorem 2.1 of \cite{PWW17} and Theorem 1 of \cite{DDW23} under the condition $\ \sup_{k\ge 1}\{n_{k}\}<+\infty$, it also generalizes Theorem 1(2) of \cite{YD20} and  Corollary 1 of \cite{LLMX}.
\end{remark}
\bigskip

\begin{thm}\label{thm2}
Let $E\in \mathcal{M}(I, \{n_{k}\}, \{c_{k}\})$, if $\ \sup_{k\ge 1}\{n_{k}\}<+\infty$ , then
\begin{equation*}
\dim_{L}E\le\liminf\limits_{l\rightarrow +\infty}\inf\limits_{k\geq 1}\frac{\log n_{k+1}\cdots n_{k+l}}
{-\log c_{k+1}\cdots c_{k+l}}.
\end{equation*}
\end{thm}
\bigskip

\begin{thm}\label{thm3}
Let $E\in\mathcal{C}(J,  \{c_{k}\},\{n_{k}\},\{a_{k}\})$, then
\begin{equation*}
\dim_{A}^{\theta}E=\limsup\limits_{k\rightarrow +\infty}\frac{\log n_{k+1}\cdots n_{l(\theta,k)}}
{(1-\frac{1}{\theta})\log \delta_{k}},
\end{equation*}
where $l(\theta,k)=\max \{l\in \mathbb{Z}^{+}: \delta_{l}\geq (\delta_{k})^{\frac{1}{\theta}}\}$ \ for any \ $k \geq 1$.
\end{thm}
\bigskip

\begin{thm}\label{thm4}
Let $E\in\mathcal{C}(J,  \{c_{k}\},\{n_{k}\},\{a_{k}\})$, then
\begin{equation*}
\dim_{L}^{\theta}E=\liminf\limits_{k\rightarrow +\infty}\frac{\log n_{k+1}\cdots n_{l(\theta,k)}}
{(1-\frac{1}{\theta})\log \delta_{k}},
\end{equation*}
where $l(\theta,k)=\max \{l\in \mathbb{Z}^{+}: \delta_{l}\geq (\delta_{k})^{\frac{1}{\theta}}\}$ \ for any \ $k \geq 1$.
\end{thm}

\begin{remark}
Any homogeneous perfect set $E \in \mathcal{J}(I, \{n_{k}\}, \{c_{k}\}, \{\eta_{k,j}\} )$ with $c_{\ast}=\inf_{k}\{c_{k}\}>0$ is a Cantor-like set,  where $E\in\mathcal{C}(I,  \{c_{k}\},\{n_{k}\},\{a_{k}\})$ with $a_{k}=0$ for all $k\ge 1$. Thus Theorem \ref{thm3} generalizes Theorem 2 of \cite{YD20} and Theorem \ref{thm4} generalizes Theorem 2 of \cite{YL22}.
\end{remark}

\bigskip
\bigskip

\section{Assouad dimension and lower dimension of homogeneous Moran sets}
To prove Theorem \ref{thm1} and Theorem \ref{thm2}, we need the equivalent definitions of the Assouad dimension and the lower dimension as follows, see in \cite{C19}.
\begin{lemma}{\rm(equivalent definitions \cite{C19})} \label{eq}

Let $F$ be a compact set in a metric space $(X, d)$, then
\begin{equation}
\begin{split}
\dim_{A}F=\inf \bigg\{s\geq 0:\ & \text{there exist} \ b,\ c>0 \ \text{such that for any}\ \ 0<r<b, ~~~~\text{any}\ \ r<R< |F|, \\ &\text{and any}\ x\in F, ~~~~~ N_{r}\big(B(x,R)\cap F\big)\leq c\big(\frac{R}{r}\big)^{s} \bigg\},   \label{Assouad11}
\end{split}
\end{equation}
\begin{equation}
\begin{split}
\dim_{L}F=\sup \bigg\{s\geq 0:\ & \text{there exist} \ b,\ c>0 \ \text{such that for any}\ \ 0<r<b, \text{any}\ \ r<R< |F|, \\ &\text{and any}\ x\in F, ~~~~~ N_{r}\big(B(x,R)\cap F\big)\geq c\big(\frac{R}{r}\big)^{s} \bigg\}.   \label{Assouad21}
\end{split}
\end{equation}
\end{lemma}
\bigskip

We start to study the Assouad dimensions and the lower dimensions of homogeneous Moran sets under the condition $\sup_{k\ge 1}\{n_{k}\}<+\infty$.  Let $M=\sup_{k\ge 1}\{n_{k}\}$.
\subsection{Proof of Theorem \ref{thm1}}
To prove Theorem \ref{thm1}. First, we proof $$\dim_{A}E\ge \limsup\limits_{l\rightarrow +\infty}\sup\limits_{k\geq 1}\frac{\log n_{k+1}\cdots n_{k+l}}
{-\log c_{k+1}\cdots c_{k+l}}.$$

 If $s>\dim_{A}E$, by (\ref{Assouad11}),  there  exist constants $ b, C>0$, such that for any $0<r< b$, any $r<R< |E|$, and any $x\in E$,
$$
N_{r}\big(B(x,R)\cap E\big)\leq C\big(\frac{R}{r}\big)^{s}. $$
Then there exist $L_{1}>0$ and $C_{1}>0$, such that for any $k\ge 1$, $l\ge L_{1}$, $R=\delta_{k}$ and $r=\delta_{k+l}$,
\begin{equation}\label{1}
N_{r}\big(B(x,R)\cap E\big)\leq C_{1}\big(\frac{R}{r}\big)^{s}=C_{1}(c_{k+1}\cdot\cdot\cdot c_{k+l})^{-s}
\end{equation}
for any $x\in E$.

Notice that for any $x\in E$, $B(x,R)$ contains at least one $(k+1)-$level basic interval of $E$, and $B(x,r)$ meets at most four $(k+l)-$level basic intervals of $E$. Thus
\begin{equation}\label{2}
\begin{split}
N_{r}\big(B(x,R)\cap E\big)&\geq \frac{n_{k+2}\cdot\cdot\cdot n_{k+l}}{4}\\&\geq \frac{1}{4M}(n_{k+1}\cdot\cdot\cdot n_{k+l})
\end{split}
\end{equation}
for any $x\in E$.

Combine (\ref{1}) and (\ref{2}), we have
\begin{equation}
\frac{1}{4M}(n_{k+1}\cdot\cdot\cdot n_{k+l})\leq C_{1}(c_{k+1}\cdot\cdot\cdot c_{k+l})^{-s}
\end{equation}
for any $k\ge 1$ and $l\ge L_{1}$, which implies that there is a constant $C_{2}$, such that
$$ \frac{\log{n_{k+1}\cdot\cdot\cdot n_{k+l}}}{-\log{c_{k+1}\cdot\cdot\cdot c_{k+l}}}\leq s+ \frac{C_{2}}{\log{c_{k+1}\cdot\cdot\cdot c_{k+l}}}$$
for any $k\ge 1$ and $l\ge L_{1}$.

Let $l\rightarrow +\infty$, we have $$\limsup\limits_{l\rightarrow +\infty}\sup\limits_{k\geq 1}\frac{\log n_{k+1}\cdots n_{k+l}}
{-\log c_{k+1}\cdots c_{k+l}}\leq s.$$ By the arbitrariness of $s>\dim_{A}E$, we have $$\dim_{A}E\ge\limsup\limits_{l\rightarrow +\infty}\sup\limits_{k\geq 1}\frac{\log n_{k+1}\cdots n_{k+l}}
{-\log c_{k+1}\cdots c_{k+l}}.$$

Now we proof $$\dim_{A}E\le \limsup\limits_{l\rightarrow +\infty}\sup\limits_{k\geq 1}\frac{\log n_{k+1}\cdots n_{k+l}}
{-\log c_{k+1}\cdots c_{k+l}}.$$

For any $s>\limsup\limits_{l\rightarrow +\infty}\sup\limits_{k\geq 1}\frac{\log n_{k+1}\cdots n_{k+l}}
{-\log c_{k+1}\cdots c_{k+l}}$, there is a positive integer $L_{2}$, such that for any $l\ge L_{2}$ and $k \ge 1$, we have
$$
s>\frac{\log{n_{k+1}\cdot\cdot\cdot n_{k+l}}}{-\log{c_{k+1}\cdot\cdot\cdot c_{k+l}}},
$$
which implies that
\begin{equation}\label{3}
n_{k+1}\cdot\cdot\cdot n_{k+l}\le (c_{k+1}\cdot\cdot\cdot c_{k+l})^{-s}.
\end{equation}

For any $0<R<\delta_{1}$, there exists a positive integer $k\geq 2$, such that
\begin{equation}\label{4}
\delta_{k}<R\le\delta_{k-1}.
\end{equation}

For any $0<r<\delta_{k+L_{2}}$, there exists a positive integer $l \geq L_{2}$, such that
\begin{equation}\label{5}
\delta_{k+l+1}< r\le\delta_{k+l}.
\end{equation}
 Then for any $x\in E$, $B(x,R)$ meets at most four $(k-1)-$level basic intervals of $E$, and $B(x,r)$ contains at least one $(k+l+1)-$level basic interval of $E$. Notice that $l \geq L_{2}$, by (\ref{3}), (\ref{4}) and (\ref{5}), we have
\begin{equation}\label{6}
\begin{split}
N_{r}\big(B(x,R)\cap E\big)&\leq 4n_{k}\cdot\cdot\cdot n_{k+l+1}\\&\leq 4M^{2}(n_{k+1}\cdot\cdot\cdot n_{k+l})\\&\leq 4M^{2}(c_{k+1}\cdot\cdot\cdot c_{k+l})^{-s}\\&\leq4M^{2}(\frac{R}{r})^{s}
\end{split}
\end{equation}
for any $x\in E$.

For any $\delta_{k+L_{2}}\le r<R<\delta_{1}$, $B(x,r)$ contains at least one $(k+L_{2}+1)-$level basic interval of $E$. Then
\begin{equation}\label{7}
\begin{split}
N_{r}\big(B(x,R)\cap E\big)&\leq 4n_{k}\cdots n_{k+L_{2}+1}\\
&\leq 4M^{L_{2}+2}\\
&\leq 4M^{L_{2}+2}\big(\frac{R}{r}\big)^{s}
\end{split}
\end{equation}
for any $x\in E$.

Combine (\ref{6}) and (\ref{7}), we have $\dim_{A}E\leq s$. By the arbitrariness of $s>\limsup\limits_{l\rightarrow +\infty}\sup\limits_{k\geq 1}\frac{\log n_{k+1}\cdots n_{k+l}}
{-\log c_{k+1}\cdots c_{k+l}}$, we obtain  $$\dim_{A}E\le \limsup\limits_{l\rightarrow +\infty}\sup\limits_{k\geq 1}\frac{\log n_{k+1}\cdots n_{k+l}}
{-\log c_{k+1}\cdots c_{k+l}}$$ and complete the proof of Theorem \ref{thm1}.
\bigskip

\subsection{Proof of Theorem \ref{thm2}}
 If $s<\dim_{L}E$, by (\ref{Assouad21}),  there  exist constants $ b, C>0$, such that for any $0<r< b$, any $r<R< |E|$, and any $x\in E$,

$$
N_{r}\big(B(x,R)\cap E\big)\geq C\big(\frac{R}{r}\big)^{s}. $$
Then there exist  $L_{3}>0$ and $C_{3}>0$, such that for any $k\ge 1$, $l\ge L_{3}$, $R=\delta_{k}$ and $r=\delta_{k+l}$,
\begin{equation}\label{8}
N_{r}\big(B(x,R)\cap E\big)\geq C_{3}\big(\frac{R}{r}\big)^{s}=C_{3}(c_{k+1}\cdot\cdot\cdot c_{k+l})^{-s}
\end{equation}
for any $x\in E$.

Notice that for any $x\in E$, $B(x,R)$ meets at most four $k-$level basic intervals of $E$, and $B(x,r)$ contains at least one $(k+l+1)-$level basic interval of $E$. Thus
\begin{equation}\label{9}
\begin{split}
N_{r}\big(B(x,R)\cap E\big)&\leq 4n_{k+1}\cdot\cdot\cdot n_{k+l+1}\\&\leq 4M(n_{k+1}\cdot\cdot\cdot n_{k+l})
\end{split}
\end{equation}
for any $x\in E$.

Combine (\ref{8}) and (\ref{9}), we have
\begin{equation}
4M(n_{k+1}\cdot\cdot\cdot n_{k+l})\geq C_{3}(c_{k+1}\cdot\cdot\cdot c_{k+l})^{-s}
\end{equation}
for any $k\ge 1$ and $l\ge L_{3}$, which implies that there is a constant $C_{4}$, such that
$$ \frac{\log{n_{k+1}\cdot\cdot\cdot n_{k+l}}}{-\log{c_{k+1}\cdot\cdot\cdot c_{k+l}}}\geq s+ \frac{C_{4}}{\log{c_{k+1}\cdot\cdot\cdot c_{k+l}}}$$
for any $k\ge 1$ and $l\ge L_{3}$.

Let $l\rightarrow +\infty$, we have $$\liminf\limits_{l\rightarrow +\infty}\inf\limits_{k\geq 1}\frac{\log n_{k+1}\cdots n_{k+l}}
{-\log c_{k+1}\cdots c_{k+l}}\geq s.$$ By the arbitrariness of $s<\dim_{L}E$, we have $$\dim_{L}E\le\liminf\limits_{l\rightarrow +\infty}\inf\limits_{k\geq 1}\frac{\log n_{k+1}\cdots n_{k+l}}
{-\log c_{k+1}\cdots c_{k+l}}$$
and complete the proof of Theorem \ref{thm2}.
\bigskip

\begin{remark}
To  estimate the lower dimention, we cannot use the dual method of the estimatation of the Assouad
dimension to estimate the lower bound since $c_{\ast}=\inf_{k}\{c_{k}\}=0$ implies that $\frac{R}{r}$ has not uniformly  upper bound when $\delta_{k+L_{2}}\le r<R$. Maybe $\liminf\limits_{l\rightarrow +\infty}\inf\limits_{k\geq 1}\frac{\log n_{k+1}\cdots n_{k+l}}
{-\log c_{k+1}\cdots c_{k+l}}$ is not the supremum of the lower dimension
of a homogeneous Moran set $E\in \mathcal{M}(I, \{n_{k}\}, \{c_{k}\})$ with $\ \sup_{k\ge 1}\{n_{k}\}<+\infty$, but we have not find examples to support it.
\end{remark}
\bigskip
\bigskip

\section{Assouad spectrum and lower spectrum of Cantor-like sets}
Let $X$ be a metric space and  $s\ge 0$, $A$ be a compact set in $X$, if for any $x\in A$ and $0<r\leq |A|$, there exist $c \geq 1$  and a Borel measure $\mu$ supported on $A$ satisfying
\begin{equation*}
c^{-1}r^{s}\leq \mu\big(B(x,r)\big)\leq cr^{s},
\end{equation*}
then we call $A$ an Ahlfors $s-$regular set\cite{M95}.

The Ahlfors regular sets  contain many typical fractal sets, such as the self-similar and the self conformal sets which satisfy the open set condition,  they play effective roles in the study of the fractal geometry. For any Ahlfors $s-$regular set $F$, the Assouad dimension of $F$ is $s$, which is coincided with the Hausdorff dimension  and  the packing dimension of $F$, see in  \cite{Luu}.

 Reference \cite{M95} showed that for any   Ahlfor $s-$regular set $A$, there are constants $\lambda \geq 1$,\ $a\in (0,1)$, and $1<\alpha \leq \beta <+\infty$ which satisfy that for any $x,\ x_{1},\ x_{2} \in A$ and $r\in (0, |A|)$,
\begin{equation}
\lambda^{-1}\leq \frac{\mu\big(B(x_{1},r)\big)}{\mu\big(B(x_{2},r)\big)}\leq \lambda, \label{homogeneity}
\end{equation}
\begin{equation}
\alpha \leq \frac{\mu\big(B(x,r)\big)}{\mu\big(B(x,ar)\big)}\leq \beta \label{doubling},
\end{equation}
which implies that $\mu$ is a doubling measure on the set $A$.
For any set $E$  satisfying (\ref{homogeneity}) and (\ref{doubling}), \ L\"{u}, Lou, Wen and Xi \cite{LLWX} showed that $\dim_{H}E=\liminf\limits_{r\rightarrow 0} h(r), \ \ \dim_{P}E=\limsup\limits_{r\rightarrow 0} h(r)$ , where $h(r):\ (0, \epsilon)\rightarrow \mathbb{R}^{+}$, which is called the scale function of $E$, is the function such that
\begin{equation}\label{sc1}
0<\inf\limits_{r<\epsilon}h(r)\leq \sup\limits_{r<\epsilon}h(r)<+\infty
\end{equation}
for a positive constant $\epsilon$. Futhermore,  for any $x\in E$ and $ 0<r<\epsilon$, there exists a constant $c\ge 0$ with
\begin{equation}\label{sc2}
\big|h(r)\log r-\log \mu\big(B(x,r)\big)\big|\leq c.
\end{equation}
 For an  Ahlfors  $s-$regular set,  we can take $h(r)=s$.

\bigskip
\bigskip
We are ready to  prove of Theorem \ref{thm3} and Theorem \ref{thm4}. The next propositions give the Assouad spectrum and the lower spectrum by the scale function for the sets satisfying (\ref{homogeneity})  and (\ref{doubling}).
\begin{proposition}\label{pro1}{\rm \cite{YD20}}
Let $F$ be a set in a metric space $X$, and $F$ satisfies (\ref{homogeneity}) and (\ref{doubling}) with a scale function $h(r)$, then
\begin{equation}
\dim_{A}^{\theta}F=\limsup\limits_{r\rightarrow 0} \bigg|\frac{h(r)\log r-h(r^{\frac{1}{\theta}})\log r^{\frac{1}{\theta}}}{(1-\frac{1}{\theta})\log r} \bigg|.
\end{equation}
\bigskip

\end{proposition}
\begin{proposition}\label{pro2}{\rm \cite{YL22}}
Let $F$ be a set in a metric space $X$, and $F$ satisfies (\ref{homogeneity}) and (\ref{doubling}) with a scale function $h(r)$, then
\begin{equation}
\dim_{L}^{\theta}F=\liminf\limits_{r\rightarrow 0} \bigg|\frac{h(r)\log r-h(r^{\frac{1}{\theta}})\log r^{\frac{1}{\theta}}}{(1-\frac{1}{\theta})\log r} \bigg|.
\end{equation}
\end{proposition}
\bigskip

The next proposition  shows a  scale function of a Cantor-like set.
\begin{proposition}\label{pro3}{\rm \cite{LLMX}}
Let $E\in \mathcal{C}(J, \{n_{k}\}, \{c_{k}\}, \{a_{k}\})$ , then $F$ satisfies (\ref{homogeneity}) and (\ref{doubling}) with
a scale function
\begin{equation*}
h(r)=\frac{\log N_{k}}{-\log \delta_{k}}
\end{equation*}
for any $k\ge 1$ and $\delta_{k}|J|<r\leq \delta_{k-1}|J|$.
\end{proposition}

\bigskip

\bigskip
Now we finish the proof of Theorem \ref{thm3} and Theorem \ref{thm4}.

Let $E\in \mathcal{C}(J, \{n_{k}\}, \{c_{k}\}, \{a_{k}\})$. For any $r\in (0,\min\{1,|J|\})$, there is $k\geq 1$ satisfying
\begin{equation*}
\delta_{k}|J|<r\leq \delta_{k-1}|J|.
\end{equation*}
 Since $l(\theta,k)=\max \{l\in \mathbb{Z}^{+}: \delta_{l}\geq (\delta_{k})^{\frac{1}{\theta}}\}$, we have $l(\theta,k)\geq k$
and
\begin{equation*}
\delta_{l(\theta,k)+1}<(\delta_{k})^{\frac{1}{\theta}}\leq \delta_{l(\theta,k)}.
\end{equation*}
Let $l^{\ast}(\theta,r)$ be the non-negative integer satisfying \begin{equation*}
\delta_{l^{\ast}(\theta,r)+1}<(\frac{r}{|J|})^{\frac{1}{\theta}}\leq \delta_{l^{\ast}(\theta,r)}.
\end{equation*}
Notice that $$\delta_{l(\theta,k)+1}<(\delta_{k})^{\frac{1}{\theta}}<(\frac{r}{|J|})^{\frac{1}{\theta}} $$
and
$$(\frac{r}{|J|})^{\frac{1}{\theta}}\le (\delta_{k-1})^{\frac{1}{\theta}}\le \delta_{l(\theta,(k-1))},$$
then we have $$l(\theta,(k-1))\le l^{\ast}(\theta,r)\le l(\theta, k).$$

Let $l^{\ast\ast}(\theta,r)$ be the non-negative integer satisfying \begin{equation*}
\delta_{l^{\ast\ast}(\theta,r)+1}|J|<r^{\frac{1}{\theta}}\leq \delta_{l^{\ast\ast}(\theta,r)}|J|.
\end{equation*}
Since $0<c_{\ast}<c^{\ast}=\sup_{k}\{c_{k}\}<1$ and  $\sum_{k=1}^{+\infty}
a_{k}<+\infty$, we have $\sup_{k\ge 1}\{n_{k}\}<+\infty$, $\sup_{k\ge1}\{l(\theta, k)-l(\theta,(k-1))\}<+\infty$ and $\sup_{0<r<\min\{1,|J|\}}\{|l^{\ast\ast}(\theta,r)-l^{\ast}(\theta,r)|\}<+\infty$, which together with Proposition \ref{pro1}, Proposition \ref{pro2} and Proposition \ref{pro3} yields
\begin{equation*}
\begin{split}
\dim_{A}^{\theta}E&=\limsup\limits_{r\rightarrow 0} \bigg|\frac{h(r)\log r-h(r^{\frac{1}{\theta}})\log r^{\frac{1}{\theta}}}{(1-\frac{1}{\theta})\log r} \bigg|\\&=\limsup\limits_{r\rightarrow 0} \bigg|\frac{h(r)\log r-\frac{\log N_{l^{\ast\ast}(\theta,r)+1}}{-\log \delta_{l^{\ast\ast}(\theta,r)+1}}\log \delta_{l^{\ast\ast}(\theta,r)+1}}{(1-\frac{1}{\theta})\log r} \bigg|\\&=\limsup\limits_{r\rightarrow 0} \bigg|\frac{h(r)\log r-\frac{\log N_{l^{\ast}(\theta,r)+1}}{-\log \delta_{l^{\ast}(\theta,r)+1}}\log \delta_{l^{\ast}(\theta,r)+1}}{(1-\frac{1}{\theta})\log r} \bigg|\\&=\limsup\limits_{k\rightarrow +\infty} \bigg|\frac{\frac{\log N_{k}}{-\log \delta_{k}}\log \delta_{k}-\frac{\log N_{l(\theta,k)}}{-\log \delta_{l(\theta,k)}}\log \delta_{l(\theta,k)}}{(1-\frac{1}{\theta})\log \delta_{k}} \bigg|\\&=\limsup\limits_{k\rightarrow +\infty}\frac{\log n_{k+1}\cdots n_{l(\theta,k)}}
{(1-\frac{1}{\theta})\log \delta_{k}}
\end{split}
\end{equation*}
and

\begin{equation*}
\begin{split}
\dim_{L}^{\theta}E&=\liminf\limits_{r\rightarrow 0} \bigg|\frac{h(r)\log r-h(r^{\frac{1}{\theta}})\log r^{\frac{1}{\theta}}}{(1-\frac{1}{\theta})\log r} \bigg|\\&=\liminf\limits_{r\rightarrow 0} \bigg|\frac{h(r)\log r-\frac{\log N_{l^{\ast\ast}(\theta,r)+1}}{-\log \delta_{l^{\ast\ast}(\theta,r)+1}}\log \delta_{l^{\ast\ast}(\theta,r)+1}}{(1-\frac{1}{\theta})\log r} \bigg|\\&=\liminf\limits_{r\rightarrow 0} \bigg|\frac{h(r)\log r-\frac{\log N_{l^{\ast}(\theta,r)+1}}{-\log \delta_{l^{\ast}(\theta,r)+1}}\log \delta_{l^{\ast}(\theta,r)+1}}{(1-\frac{1}{\theta})\log r} \bigg|\\&=\liminf\limits_{k\rightarrow +\infty} \bigg|\frac{\frac{\log N_{k}}{-\log \delta_{k}}\log \delta_{k}-\frac{\log N_{l(\theta,k)}}{-\log \delta_{l(\theta,k)}}\log \delta_{l(\theta,k)}}{(1-\frac{1}{\theta})\log \delta_{k}} \bigg|\\&=\liminf\limits_{k\rightarrow +\infty}\frac{\log n_{k+1}\cdots n_{l(\theta,k)}}
{(1-\frac{1}{\theta})\log \delta_{k}}.
\end{split}
\end{equation*}
\bigskip

\textbf{Acknowledgements} The authors thank  the reviewers for their  helpful comments and suggestions.


\bigskip
\bigskip
\begin{thebibliography}{99}


\bibitem{As1} P.\ Assouad,
\newblock Espaces m\'{e}triques, plongements, facteurs, Th\`{e}se de doctorat d'Etat,
\newblock Publ. Math. Orsay 223--7769, Univ. Paris XI, Orsay (1977).

\bibitem{As2} P.\ Assouad,
\newblock \'{E}tude d'une dimension m\'{e}trique li\'{e}e \`{a} la possibilit\'{e} de plongements dans $\mathbb{R}^{n}$,
\newblock C. R. Acad. Sci. Paris S\'{e}r. A-B  288 (1979)  731--734.
\bibitem {Fra21} J.\ Fraser,
\newblock Assouad dimension and fractal geometry,
\newblock Cambridge University Press, 2021.

\bibitem{Fra14} J.\  Fraser,
\newblock Assouad type dimensions and homogeneity of fractals,
\newblock Trans. Amer. Math. Soc. 366 (2014) 6687--6733.

\bibitem{FT18} J.\  Fraser, M.\ Todd,
\newblock Quantifying inhomogeneity in fractal sets,
\newblock Nonlinearity 31(4) (2018) 1313--1330.

\bibitem{Luu} J.\ Luukkainen,
\newblock Assouad dimension: antifractal metrization, porous sets, and homogeneous measures,
\newblock J.\ Korean Math.\ Soc. 35(1998) 23--76.
\bibitem{Ma11} J.\ Mackay,
\newblock Assouad dimension of self-affine carpets,
\newblock Conform Geom. Dyn.  15(2011) 177--187.

\bibitem{LLMX} W.\ W.\ Li, W.\ X.\ Li, J.\ J.\ Miao, L.\ F.\ Xi,
\newblock On the Assouad dimension of Moran sets and Cantor-like sets,
\newblock Front.\ Math.\ China 11(2016) 705--722.
\bibitem{C19}  H.\ P.\ Chen, \newblock Assouad dimensions and spectra of Moran cut-out sets, \newblock Chaos. \  Solution.\  Fract 199 (2019) 310--317.


\bibitem{PWW17} F.\ J.\ Peng, W.\ Wang, S.\ Y.\ Wen,
\newblock On Assouad dimension of products,
\newblock Chaos. \  Solution.\ Fract 104 (2017) 192--197.
\bibitem{YD20} J.\ J.\ Yang, Y.\ L.\ Du,
\newblock Assouad dimension and spectrum of homogeneous perfect sets,
\newblock Fractals  28(7) (2020)  2050132.
\bibitem{DDW23} Y.\ X.\ Dai, J.\ M.\ Dong, C.\ Wei,
\newblock  The continuity of dimensions and quasisymmetrical equivalence
of parameterized homogeneous Moran sets,
\newblock J.\ Math.\ Anal.\ Appl.  518 (2023), 126783.
\bibitem {Lar1} D.\ G.\ Larman,
\newblock A new theory of dimension.
\newblock Proc.\ Lond.\ Math.\ Soc. 17 (1967), 178--192.




\bibitem{CWW17}  H.\ P.\ Chen, M.\ Wu, C. Wei. \newblock  Lower dimensions of some fractal sets, \newblock J.\ Math.\ Anal.\ Appl. 445(2) (2017) 1022--1036.
\bibitem{YL22} J.\ J.\ Yang, Y.\ Z.\ Li, R.\ Hu,
\newblock Lower dimension and spectrum of homogeneous perfect sets,
\newblock Fractals  30(9) (2022)  2250184.

\bibitem{FH18} J.\ Fraser, Y.\ Han,
\newblock New dimension spectra: finer information on scaling and homogeneity,
\newblock  Adv. Math. 329 (2018) 273--328.

\bibitem{CWC} H.\ P.\ Chen, M.\ Wu, Y.\ Y.\ Chang,
\newblock Lower Assouad type dimensions of uniformly perfect sets in doubling metric spaces.
\newblock Fractals 28(2) (2020)  2050039.
\bibitem{Wen00} Z.\ Y.\ Wen,
\newblock Mathematical Foundation of Fractal Geometry,
\newblock Shanghai Scientific Technological Education Publishing House, 2000.


\bibitem{FWW97} D.\ J.\ Feng,  Z.\ Y.\ Wen, J.\ Wu,
\newblock Some dimensional results for homogeneous Moran sets,
\newblock Sci. \ China Ser.\ A. 40 (1997) 475--482.
\bibitem{KLS} Y.\ Z.\ Li, X.\ H.\ Fu, J.\ J.\ Yang,
\newblock Quasisymmetrically packing-minimal perfect Moran sets,
\newblock Fractals 29(2) (2021)  2150043.



\bibitem{LW11} J.\ J.\ Li, M.\ Wu,
\newblock Pointwise dimensions of general Moran measures with open set condition,
\newblock Sci. China Math. 54(4) (2011) 699--710.



\bibitem{M95} P.\ Mattila,
\newblock Geometry of Sets and Measures in Euclidean Spaces: Fractals and Rectifiability,
\newblock Cambridge University Press, 1995.

\bibitem {LLWX} F.\ L\"{u}, M.\ L.\ Lou, Z.\ Y.\ Wen, L.\ F.\ Xi,
\newblock Bilipschitz embedding of homogeneous fractals,
\newblock J.\ Math.\ Anal.\ Appl. 432 (2015) 888--917.
\end{thebibliography}
\end{document}